\documentstyle[amssymb,12pt]{article}

\begin{document}

\newtheorem{theorem}{Theorem}{}
\newtheorem{lemma}[theorem]{Lemma}{}
\newtheorem{corollary}[theorem]{Corollary}{}
\newtheorem{conjecture}[theorem]{Conjecture}{}
\newtheorem{proposition}[theorem]{Proposition}{}
\newtheorem{axiom}{Axiom}{}
\newtheorem{remark}{Remark}{}
\newtheorem{example}{Example}{}
\newtheorem{exercise}{Exercise}{}
\newtheorem{definition}{Definition}{}

\title{Towards Drinfeld--Sokolov reduction for quantum groups}
\author{Alexey Sevostyanov
 \footnote{e-mail sevastia@mpim-bonn.mpg.de}\\ 
Max-Planck-Institute f\"{u}r mathematik, Box 7280, D-53072 Bonn,
Germany
 }
\maketitle

\begin{abstract}
In this paper we study the Poisson--Lie version of the Drinfeld--Sokolov
reduction defined in \cite{FRS}, \cite{SS}. Using the bialgebra structure
related to the new Drinfeld realization of affine quantum groups we describe
reduction in terms of constraints. This realization of reduction admits
direct quantization.

As a byproduct we obtain an explicit expression for the symplectic form
associated to the twisted Heisenberg double and calculate the moment map for
the twisted dressing action. For some class of infinite--dimensional Poisson
Lie groups we also prove an analogue of the Ginzburg--Weinstein isomorphism.
\end{abstract}

\section*{Introduction}

It is well known that the quantum Drinfeld--Sokolov reduction plays an
important role in the theory of affine Lie algebras. Let ${\frak g}$ be a
complex semisimple Lie algebra, $\widehat{{\frak g}}$ the affinization of $%
{\frak g}$. In \cite{FF} Feigin and Frenkel proved that the
Drinfeld--Sokolov reduction associated to ${\frak g}$ describes the
structure of the center of the universal enveloping algebra $U(\widehat{%
{\frak g}})$ at the critical level of the central charge.

The Poisson--Lie version of the Drinfeld--Sokolov reduction was proposed in
\cite{FRS},\cite{SS}. The quantization of this reduction is expected to play
the same role for affine quantum groups. In this paper we obtain a
realization of the Drinfeld--Sokolov reduction for Poisson--Lie groups
admitting direct quantization.

Recall that in physical terms every reduction procedure consists of two
steps: (1) imposing constraints , (2) choosing a cross--section of the
constraint surface. To quantize a system with constraints one should
quantize the underlying Poisson manifold, find quantum counterparts of the
classical constraints and apply the quantum reduction procedure to the
quantum system. All the known versions of this procedure require an explicit
description of constraints. For the Lie--Poisson version of the
Drinfeld--Sokolov reduction both the quantization of the underlying Poisson
manifold and the description of constraints are nontrivial problems.

The main observation of \cite{FRS} , \cite{SS} was that in order to perform
the Poisson--Lie version of the Drinfeld--Sokolov reduction associated to a
complex semisimple Lie algebra ${\frak g}$ one should introduce a new
quasitriangular bialgebra structure on the loop algebra ${\frak g}((z))$
arising from natural geometric considerations. The corresponding r--matrix
is obtained by adding an extra term to the standard r--matrix related to the
``new Drinfeld realization'' of affine quantum groups. This term is
essentially elliptic and may be expressed by means of theta functions (see
section \ref{srtheta}). We denote by ${\cal G}$ the corresponding Poisson
Lie group. The Poisson manifold which undergoes reduction is essentially the
dual Poisson Lie group, with its standard Poisson structure being twisted by
an automorphism satisfying certain conditions (see (\ref{compat})). In \cite
{FRS} this Poisson structure is called gauge covariant. We denote the
corresponding Poisson manifold by ${\cal G}_{p}$. This manifold is equipped
with the twisted dressing action of the Poisson--Lie group ${\cal G}$.

In the theory of Poisson group actions constraints naturally appear as
matrix coefficients of moment maps in the sense of Lu and Weinstein.
However, the twisted dressing action of the unipotent subgroup ${\cal N}%
=LN\subset {\cal G}$ used in \cite{FRS},\cite{SS} for reduction is not a
Poisson group action. In section \ref{s2} we define a map $ \mu_{\cal N}$ from 
the underlying Poisson manifold ${\cal G}_p$ to the opposite unipotent group 
$\overline{\cal N}$ which forms a dual pair together with the canonical projection onto 
the quotient ${\cal G}_p / {\cal N}$ and serves as a substitution of the moment map.

As a matter of fact,  it is not really necessary to quantize the new
elliptic bialgebra structure. It is well known that at least in the finite
dimensional case quantizations of different bialgebra structures are
isomorphic as algebras (see \cite{D}). The classical counterpart of this
statement for Lie--Poisson groups is called the Ginzburg--Weinstein
isomorphism \cite{GW}. We show that the same statement holds for some class
of infinite--dimensional Poisson--Lie groups. Applying a simple form of the
Ginzburg--Weinstein isomorphism found in \cite{A}, we prove that the gauge
covariant Poisson structures corresponding to the new bialgebra structure
and to the one related to Drinfeld's ``new realization'' of affine quantum
groups are isomorphic. This allows us to use the latter for reduction.
Surprisingly, the description of constraints in this realization does not
contain elliptic functions. Quantum counterparts of these constraints have
been defined in \cite{S1}.

Remarkably, the quantum constraints are of the first class, i.e. they form a
subalgebra that possesses a character fixing the values of the constraints.
An appropriate reduction technique for such systems of constraints has been
developed by the author in \cite{S2}.

\subsubsection*{Acknowledgments}

The author would like to thank M.A. Semenov--Tian--Shansky for useful
discussions and A. Alekseev for valuable advice in symplectic geometry. I am
also grateful to M. Golenishcheva-Kutuzova for careful reading of the text.

\section{The Poisson geometry of the twisted dressing action}

\label{s1}

In this section we develop the moment map technique for the twisted dressing
action. The moment map is important both to describe the constraints and to
prove the existence of the Ginzburg--Weinstein isomorphism (see \cite{A}).

Let $G\times M\rightarrow M$ be a Poisson group action of a Poisson Lie
group $G$ on a Poisson manifold $M$ possessing a moment map $\mu $. Then $%
\mu $ maps the manifold $M$ into the dual Poisson--Lie group $G^{*}$. The
most important particular example of a Poisson group action is the dressing
action of a Poisson Lie group $G$ on the dual group $G^{*}$; in that case
the moment map is the identity mapping. For quasitriangular Poisson Lie
groups the easiest way to obtain this action is to consider $G^{*}$ as the
reduced Poisson manifold for the Heisenberg double $D$ of $G$ with respect
to the left Poisson group action of $G$. Recall that the Heisenberg double
is isomorphic to $G\times G$ as a manifold and $G$ acts on the product $%
G\times G$ by left and right diagonal translations. The right Poisson group
action of $G$ on the Heisenberg double generates the dressing action on $%
G^{*}$.

For the needs of the Drinfeld--Sokolov reduction we have to twist the
standard Poisson structure of $G^{*}$ by an automorphism satisfying certain
conditions (see (\ref{compat})). This twisted Poisson structure, called the
gauge covariant Poisson structure in \cite{FRS}, may be obtained by a
reduction from the twisted Heisenberg double of $G$. The twisted Heisenberg
double is equipped with left and right Poisson group actions of the Poisson
Lie group $G$. The Poisson reduction with respect to the left action of $G$
yields a gauge covariant Poisson structure on $G^{*}$. Then the right
Poisson group action of the same group generates the twisted dressing action
of $G$ on $G^{*}$ called the gauge action in \cite{FRS}.

The Poisson structure on the twisted Heisenberg double is nondegenerate on
an open dense subset . Following \cite{AM}, we obtain an explicit expression
for the corresponding symplectic form. Using this formula we calculate the
moment maps for the right and left Poisson group actions of $G$ on the
twisted Heisenberg double. In contrast with the untwisted case, these maps
are neither Poisson nor equivariant. As a consequence, we get the moment map
for the gauge action of $G$ on $G^{*}$. This result will be applied to the
Drinfeld--Sokolov reduction in the next section.

\subsection{Factorizable Lie bialgebras and their doubles}

Let us recall some basic notions on Poisson Lie groups (see \cite{D}, \cite
{fact} , \cite{dual}). Let $G$ be a Lie group equipped with a Poisson
bracket , ${\frak g}$ its Lie algebra. $G$ is called a Poisson Lie group if
the multiplication $G\times G\rightarrow G$ is a Poisson map. A Poisson
bracket satisfying this axiom is degenerate and , in particular, is
identically zero at the unit element of the group. Linearizing this bracket
at the unit element we get the structure of a Lie algebra in the space $%
T_{e}^{*}G\simeq {\frak g}^{*}$. The pair (${\frak g},{\frak g}^{*})$ is
called the tangent bialgebra of $G$. (${\frak g},{\frak g}^{*})$ is called a
{\em factorizable Lie bialgebra }if the following conditions are satisfied
\cite{fact} , \cite{D}:

\begin{enumerate}
\item  ${\frak g}${\em \ is equipped with a fixed nondegenerate invariant
scalar product} $\left\langle \cdot ,\cdot \right\rangle $.

We shall always identify ${\frak g}^{*}$ and ${\frak g}$ by means of this
scalar product.

\item  {\em The dual Lie bracket on }${\frak g}^{*}\simeq {\frak g}${\em \
is given by}
\begin{equation}
\left[ X,Y\right] _{*}=\frac{1}{2}\left( \left[ rX,Y\right] +\left[
X,rY\right] \right) ,X,Y\in {\frak g},  \label{rbr}
\end{equation}

{\em where }$r\in End\ {\frak g}${\em \ is a skew symmetric linear operator
(classical r-matrix).}

\item  $r${\em \ satisfies} {\em the} {\em modified classical Yang-Baxter
identity:}
\begin{equation}
\left[ rX,rY\right] -r\left( \left[ rX,Y\right] +\left[ X,rY\right] \right)
=-\left[ X,Y\right] ,\;X,Y\in {\frak g}{\bf .}  \label{cybe}
\end{equation}
\end{enumerate}

Define operators $r_\pm \in End\ {\frak g}$ by
\[
r_{\pm }=\frac 12\left( r\pm id\right) .
\]

Then the classical Yang--Baxter equation implies that $r_{\pm }$ , regarded
as a mapping from ${\frak g}^{*}$ into ${\frak g}$ , is a Lie algebra
homomorphism. Moreover, $r_{+}^{*}=-r_{-},$\ and $r_{+}-r_{-}=id.$

The double of a factorizable Lie bialgebra admits the following explicit
description (cf.. \cite{dual}, \S 2) . Put ${\frak {d}}={\frak g\oplus {g}}$
(direct sum of two copies).The mappings 
\begin{eqnarray}
{\frak {g}}^{*} &\rightarrow &{\frak {d}}~~~:X\mapsto
(X_{+},~X_{-}),~~~X_{\pm }~=~r_{\pm }X,  \label{imbd} \\
{\frak {g}} &\rightarrow &{\frak {d}}~~~:X\mapsto (X,~X)  \nonumber
\end{eqnarray}
are Lie algebra embeddings. Thus we may identify ${\frak g^{*}}$ and ${\frak %
g}$ with Lie subalgebras in ${\frak {d}}$. Equip ${\frak {d}}$ with the
scalar product

\begin{equation}\label{prsq}
\langle\langle (X,X^{\prime }),(Y,Y^{\prime })\rangle\rangle = \langle
X,Y\rangle -\langle X^{\prime },Y^{\prime }\rangle .
\end{equation}

\begin{proposition}
{\bf (\cite{dual}, Proposition 2.1 )} \label{rdouble}

(i)$(%
{\frak d},{\frak g},{\frak g}^{*})$ is a Manin triple, i.e. ${\frak g}$ and ${\frak g}^{*}$ are
isotropic subalgebras with respect to the scalar product (\ref{prsq}).

(ii)${\frak d}$ is isomorphic to the double of $({\frak g},{\frak g}^{*}).$

(iii) $({\frak d},{\frak d}^{*})$ is a factorizable Lie bialgebra; the
corresponding r-matrix $r_{{\frak {d}}}\in End\,({\frak g}\oplus {\frak {g})}
$ is given by 
\begin{equation}
r_{{\frak {d}}}~~~=~~~\left( 
\begin{array}{ll}
r & -2r_{+} \\ 
2r_{-} & -r
\end{array}
\right) .
\end{equation}
\end{proposition}

The problem of classification of solutions of the classical Yang--Baxter
equation had been solved by Belavin and Drinfeld in \cite{BD} (cf. also \cite
{rmatr}). Their results may be summarized as follows.

Denote by ${\frak b}_\pm$ and ${\frak n}_\mp$ the image and the kernel of
the operator $r_\pm $: ${\frak b}_\pm = Im~r_\pm,{\frak n}_\mp = Ker~r_\pm $

\begin{theorem}
{\bf (Belavin--Drinfeld)}\label{rmatrtheta}

Let $({\frak g},{\frak g}^{*})$ be a factorizable Lie bialgebra. Then

(i) ${\frak b}_{\pm }\subset {\frak g}$ is a Lie subalgebra ; the subspace $%
{\frak n}_{\pm }$ is a Lie ideal in ${\frak b}_{\pm },{\frak b}_{\pm
}^{\perp }={\frak n}_{\pm }$.

(ii) The map $\theta_r :{\frak b}_{-}/{\frak n}_{-}\rightarrow {\frak b}_{+}/%
{\frak n}_{+}$ which sends the residue class of $r_{-}(X),X\in {\frak g}^{*}$%
, modulo ${\frak n}_{-}$ to that of $r_{+}(X)$ modulo ${\frak n}_{+}$ is a
well-defined isomorphism of Lie algebras. Moreover , $\theta_r $ is a unitary
operator with respect to the induced scalar product: $\theta_r \theta_r^{*}=1$.
\end{theorem}

Using part (ii) of Theorem \ref{rmatrtheta} we can describe the image of the embedding ${\frak {g}}%
^{*}\rightarrow {\frak {d}}$ as follows:

\begin{equation}  \label{g*}
{\frak g}^* = \{(X_+,X_-) \in {\frak b}_+ \oplus {\frak b}_- \subset {\frak d%
} | \overline{X}_+ = \theta_r(\overline{X}_-) \},\mbox{ where }\overline{X}%
_\pm = {X}_\pm \mbox{ mod } {\frak n}_\pm.
\end{equation}

We shall also need other properties of the subalgebras ${\frak b}_\pm$.

\begin{proposition}
{\bf (\cite{BD})} \label{bpm}

(i) ${\frak n}_{\pm }$ is an ideal in ${\frak {g}}^{*}$.

(ii) ${\frak b}_{\pm }$ is a Lie subalgebra in ${\frak {g}}^{*}$. Moreover $%
{\frak b}_{\pm }={\frak {g}}^{*}/{\frak n}_{\pm }$.

(iii) $({\frak b}_{\pm },{\frak b}_{\pm }^{*})$ is a sub-bialgebra of $(%
{\frak {g}},{\frak {g}}^{*})$ and ${\frak b}_{\pm }^{*}\simeq {\frak b}_{\mp
}$. The canonical pairing between ${\frak b}_{\mp }$ and ${\frak b}_{\pm }$
is given by

\begin{equation}
(X_{\mp },Y_{\pm })_{\pm }=\langle X_{\mp },r_{\pm }^{-1}Y_{\pm }\rangle
,~X_{\mp }\in {\frak b}_{\mp };~Y_{\pm }\in {\frak b}_{\pm }.
\end{equation}

(iv) ${\frak n}_{\pm }^{*}\simeq {\frak n}_{\mp }$ as a linear space.
\end{proposition}

\begin{remark}
\label{pmdual} If $X_{\pm }\in {\frak b}_{\pm },~(Y_{+},Y_{-})\in {\frak g}%
^{*},~Y_{\pm }=r_{\pm }Y,$ then

\begin{equation}
\langle \langle (X_{-},X_{-}),(Y_{+},Y_{-})\rangle \rangle =\langle
X_{-},Y\rangle =(X_{-},Y_{+})_{+},
\end{equation}

\begin{equation}
\langle \langle (X_{+},X_{+}),(Y_{+},Y_{-})\rangle \rangle =\langle
X_{+},Y\rangle =(X_{+},Y_{-})_{-}.
\end{equation}
\end{remark}

\subsection{Twisted Heisenberg double and its symplectic structure}

We start the study of the twisted Heisenberg double with explicit formulae
for some Poisson brackets associated with the bialgebra structure on ${\frak %
d}$.

For every Lie group $A$ with Lie algebra ${\frak a}$ we define left and
right gradients $\nabla \varphi ,\nabla ^{\prime }\varphi \in {\frak a}^{*}$%
of a function $\varphi \in C^{\infty }(A)$ by

\begin{eqnarray}  \label{grad}
\xi ( \nabla \varphi (x))\ =\left( \frac d{ds}\right) _{s=0}\varphi (e^{s\xi
}x),  \nonumber \\
\xi ( \nabla^{\prime} \varphi (x))=\left( \frac d{ds}\right) _{s=0}\varphi
(xe^{s\xi }),~~\xi \in {\frak {a}.}
\end{eqnarray}

The Poisson bracket on $D=G\times G$ associated with the bialgebra structure
on ${\frak d}$ has the form (see \cite{dual}, \S 2) :

\begin{equation}  \label{pbr}
\{ \varphi ,\psi \}=-\frac 12 \left\langle \left\langle r_{{\frak {d}}}
\nabla \varphi,\nabla \psi \right\rangle \right\rangle +\frac 12\left\langle
\left\langle r_{{\frak {d}}}\nabla^{\prime} \varphi , \nabla^{\prime} \psi
\right\rangle \right\rangle .
\end{equation}

It is well known that the bracket (\ref{pbr}) satisfies the Jacobi identity
and equips $D\ $ with the structure of a Poisson Lie group. The embeddings $%
{\frak {g}}\rightarrow {\frak {d}},~{\frak {g}}^* \rightarrow {\frak {d}}$
may be extended to homomorphisms $G\rightarrow D ,~ G^{*}\rightarrow D$. We
shall identify $G$ and $G^{*}$ with the corresponding subgroups in $D$. By
the definition of the double $G,G^{*}\subset D$ are Lie-Poisson subgroups.

Let $\sigma \in Aut~ G$ be an automorphism of $G$. We shall denote the
corresponding automorphism of ${\frak g}$ by the same letter. Assume that $%
\sigma $ satisfies the following conditions:

\begin{equation}  \label{compat}
\begin{array}{l}
1.~~ \sigma \circ r=r\circ \sigma . \\ 
\\ 
2.~~ \left\langle \sigma X,\sigma Y\right\rangle =\left\langle
X,Y\right\rangle ~~ \mbox{ for all }~~ X,Y\in {\frak g}.
\end{array}
\end{equation}

We shall associate with $\sigma $ the so called {\em twisted} Poisson
structure on $D$ (see \cite{dual}, \S 5). Let $T\in Aut\,{\frak d}$ be the
automorphism of the Lie algebra ${\frak d}$ defined by 
\[
T =~\left( \sigma ^{-1}\times id\right). 
\]
Denote by $^\sigma r_{{\frak {d}}}$ the r--matrix $r_{{\frak {d}}}$ twisted
by the automorphism $T$: 
\[
^\sigma r_{{\frak {d}}}~ =~Tr_{{\frak d}}T^{-1}~=~\left( 
\begin{array}{cc}
r & -2\sigma ^{-1}r_{+} \\ 
2r_{-}\sigma & -r
\end{array}
\right) , 
\]
and define the twisted Poisson bracket on $D$ by

\begin{equation}  \label{5.7}
\{ \varphi ,\psi \}_\sigma~~=~~\frac 12\left\langle \left\langle r_{{\frak {d%
}}} \nabla \varphi , \nabla \psi \right\rangle \right\rangle +~\frac
12\left\langle \left\langle ^\sigma r_{{\frak {d}}}\nabla^{\prime} \varphi ,
\nabla^{\prime} \psi \right\rangle \right\rangle .
\end{equation}

The Jacobi identity for (\ref{5.7}) follows from the classical Yang---Baxter
identity (\ref{cybe}) for $r_{{\frak {d}}}$.

The pair $D_{\sigma}=(D,\{,\}_\sigma )$ is called the{\em \ twisted
Heisenberg double} (for $\sigma =id$ we get the ordinary Heisenberg double).

When restricted to an open dense subset the Poisson structure of the twisted
Heisenberg double is nondegenerate. To calculate the corresponding
symplectic form we need twisted factorizations on $D$ \cite{dual}.

\begin{proposition}
\label{fact} (i) Any element $x\in {\frak d}$ admits two unique
decompositions

\begin{equation}
\begin{array}{l}
x=\eta +T\xi , \\ 
x=T\eta ^{\prime }+\xi ^{\prime }, \\ 
\eta ,\eta ^{\prime }\in {\frak g},\xi ,\xi ^{\prime }\in {\frak g}^{*}.
\end{array}
\end{equation}

(ii) In an open dense subset $D_{\sigma }^{\prime }\subset D_{\sigma }$ we
have the following factorizations :

\begin{equation}
\begin{array}{l}
d=g{g^{*}}^{T}=h^{*}h^{T}~,~{\mbox where}~d\in D_{\sigma }^{\prime }, \\ 
g,h\in G;~g^{*},h^{*}\in G^{*}, \\ 
g^{*}=(g_{+},g_{-}),h^{*}=(h_{+},h_{-}).
\end{array}
\label{factoriz}
\end{equation}
\end{proposition}

\begin{theorem}
\label{symp} $D_{\sigma }^{\prime }$ is a symplectic submanifold in $%
D_{\sigma }$. The corresponding symplectic form can be represented as
follows :

\begin{equation}
\Omega =\langle \langle \theta _{h^{*}}\otimes \theta _{g}\rangle \rangle
-\langle \langle \mu _{h}\otimes \mu _{g^{*}}\rangle \rangle ,  \label{SF}
\end{equation}

where $\theta _{h^{*}},\theta _{g}(\mu _{g^{*}},\mu _{h})$ are the universal
right--invariant (left--invariant) Maurer--Cartan forms on $G^{*}$ and $G,$
respectively; the pairing is applied to their values and lower indices
indicate group variables.
\end{theorem}

{\em Proof} of the theorem is quite similar to that for $\sigma =id$ (see 
\cite{AM}, Theorem 3).

We shall call $D_{\sigma}^{\prime}$ the principal symplectic leaf of $%
D_{\sigma}$.

One can define an action of the Drinfeld double $D$ on the Poisson manifold $%
D_{\sigma }$ which generalizes the well--known dressing action \cite{sem}.

\begin{proposition}
\label{tdact}

(i) The actions of $D$ on $D_{\sigma }$ by right and left translations are
Poisson group actions;

\begin{equation}
\begin{array}{ll}
D\times D_{\sigma }\stackrel{L}{\rightarrow }D_{\sigma }, & d^{\prime }\circ
d=d(d^{\prime }{^{T}})^{-1}; \\ 
D_{\sigma }\times D\stackrel{R}{\rightarrow }D_{\sigma }, & d\circ d^{\prime
}={d^{\prime }}^{-1}d.
\end{array}
\label{dact}
\end{equation}
where $d\in D_{\sigma },d^{\prime }\in D$.

(ii) The actions (\ref{dact}) generate Poisson group actions of the Poisson
Lie subgroups $G,G^{*}\subset D$:

\begin{equation}
G\times D_{\sigma }\stackrel{L}{\rightarrow }D_{\sigma };  \label{gactL}
\end{equation}

\begin{equation}
D_{\sigma }\times G\stackrel{R}{\rightarrow }D_{\sigma }.  \label{gactR}
\end{equation}

\begin{equation}
\begin{array}{l}
G^{*}\times D_{\sigma }\stackrel{L}{\rightarrow }D_{\sigma }; \\ 
D_{\sigma }\times G^{*}\stackrel{R}{\rightarrow }D_{\sigma }.
\end{array}
\label{g*act}
\end{equation}

(iii) The restrictions of the actions (\ref{gactL}), (\ref{gactR}), (\ref
{g*act}) to the principal symplectic leaf $D_{\sigma }^{\prime }$ possess
moment maps in the sense of Lu and Weinstein \cite{Lu}:

\begin{equation}
\begin{array}{ll}
\mu _{G}^{L}(d)=g^{*}, & \mu _{G^{*}}^{L}(d)=h, \\ 
\mu _{G}^{R}(d)=h^{*}, & \mu _{G^{*}}^{R}(d)=g,
\end{array}
\label{moment}
\end{equation}

where $g,g^{*},h,h^{*}$ are given by (\ref{factoriz}).
\end{proposition}

\begin{remark}
In general, the maps (\ref{moment}) are neither Poisson nor equivariant.
\end{remark}

{\em Proofs} of (i) and (ii) are given in \cite{sem}, Proposition 2.5.1.(See
also \cite{dual} , Proposition 5.4).

To prove (iii), let us consider , for instance , the action

\begin{equation}
D_{\sigma} \times G \stackrel{R}{\rightarrow} D_{\sigma}.
\end{equation}

Let $X\in {\frak g}$. The corresponding vector field is:

\begin{equation}
\widehat X\varphi(d) =\left({\frac{d }{ds }}\right)_{s=0}\varphi(e^{-sX}d).
\end{equation}

Therefore $\widehat{X}=(R_{d})_{*}(-X)$ , where $R_{d}$ is the operator of
right translation by $d$. According to the definition of the moment map \cite
{Lu}, we have to prove that

\begin{equation}
\Omega (\widehat X,\cdot )=-\langle \langle {\mu^R_G}^*(\theta) , X \rangle
\rangle .
\end{equation}

Formula (\ref{factoriz}) for the twisted factorization problem implies that $%
\mu _{g^{*}}(\widehat{X})=0$. Therefore, substituting $\widehat{X}$ into (%
\ref{SF}) we obtain:

\begin{equation}
\Omega (\widehat X,\cdot )=-\Omega (\cdot , \widehat X)= -\langle \langle
\theta_{h^*}\otimes_, \theta_{g}(\widehat X) \rangle \rangle .
\end{equation}

We also have $\theta _{g}(\widehat{X})=X,$ because $\widehat{X}$ is a right
invariant vector field and $G$ is a subgroup in $D$. Finally,

\begin{equation}
\Omega (\widehat X,\cdot )=-\langle \langle \theta_{h^*}\otimes_, X \rangle
\rangle .
\end{equation}

This completes the proof.

Let $B_{\pm }$ be the Lie subgroups in $G$ corresponding to the Lie
subalgebras ${\frak b}_{\pm }$. According to part (iii) of Proposition \ref
{bpm}, $B_{\pm }$ are Poisson Lie subgroups in $G$.

\begin{proposition}
\label{bact} The $G$--actions (\ref{gactL}), (\ref{gactR}) induce Poisson
group actions of the Poisson subgroups $B_{\pm }$. When restricted to the
principal symplectic leaf these actions possess moment mappings in the sense
of Lu and Weinstein:

\begin{equation}
\mu _{B_{\mp }}^{R}(d)=h_{\pm },~~\mu _{B_{\mp }}^{L}(d)=g_{\pm },
\end{equation}

where $h_{\pm },g_{\pm }$ are given by (\ref{factoriz}).
\end{proposition}

{\em Proof} follows from the previous proposition and remark \ref{pmdual}.

\subsection{Gauge covariant Poisson structures}

The gauge covariant Poisson structure used in \cite{FRS},\cite{SS} for the
Drinfeld--Sokolov reduction may be obtained from the twisted Heisenberg
double by the following construction (\cite{RIMS}, \S 3 ; \cite{dual}, \S 5)

Consider the Poisson reduction of the Poisson manifold $D_{\sigma }$ with
respect to the left action (\ref{gactL}) of the group $G$. The quotient
space $G\backslash D_{\sigma }$ may be identified with $G$, the projection
map $p:D_{\sigma }\rightarrow G$ is given by 
\[
p:(x,y)\mapsto x^{\sigma }\left( y\right) ^{-1}. 
\]

Under this identification the reduced Poisson bracket on $G$ is given by

\begin{equation}  \label{tau}
\left\{ \varphi ,\psi \right\} _\sigma =\left\langle r \nabla \varphi,\nabla
\psi \right\rangle +\left\langle r \nabla^{\prime }\varphi,\nabla^{\prime
}\psi\right\rangle -2\left\langle r_{+}^\sigma \nabla^{\prime
}\varphi,\nabla \psi\right\rangle -2\left\langle r_{-}^\sigma
\nabla\varphi,\nabla^{\prime }\psi\right\rangle ,
\end{equation}

where $r_{+}^\sigma =~\sigma \circ r_{+},~r_{-}^\sigma =r_{-}\circ \sigma
^{-1}$.

Denote by $G_\sigma =\left( G,\{,\}_\sigma \right)$ the manifold $G$
equipped with Poisson bracket (\ref{tau}).

Then the right action (\ref{gactR}) gives rise to a Poisson group action of $%
G$ on the reduced space :

\begin{equation}  \label{act}
G_\sigma \times G \rightarrow G_\sigma \;: g\circ L = (g^\sigma)^{-1}Lg.
\end{equation}

From Propositions \ref{fact} , \ref{tdact} , \ref{bact} and Theorem \ref
{symp} we deduce the following properties of the reduced Poisson manifold
and the gauge action (\ref{act}).

\begin{proposition}
\label{gprime}

(i) Elements $L\in G$ admitting a twisted factorization

\begin{equation}
L=L_{+}^{\sigma }L_{-}^{-1},(L_{+},L_{-})\in G^{*}  \label{twfact}
\end{equation}

form an open dense subset $G_{\sigma }^{\prime }$ in $G$. This factorization
is unique in a neighborhood of the unit element.

(ii) $G_{\sigma }^{\prime }$ is a Poisson submanifold in $G_{\sigma }$.

(iii) The restriction of the action (\ref{act}) to $G_{\sigma }^{\prime }$
has a moment mapping given by the identity map : $L\mapsto (L_{+},L_{-})$.

(iv) When restricted to $G_{\sigma }^{\prime }$ the induced actions of the
Poisson--Lie subgroups $B_{\mp }\subset G$ have moment maps given by :

\begin{equation}
\mu _{B_{\mp }}(L)=L_{\pm }.
\end{equation}
\end{proposition}

\section{Moment map and Drinfeld--Sokolov reduction}

\label{s2}

In this section we obtain different descriptions of the Drinfeld--Sokolov
reduction for Poisson Lie groups (see \cite{FRS} , \cite{SS} for the
definition of the reduction). Using the moment map technique developed in
the previous section, we adapt the reduction procedure for quantization.
First, we find a system of constraints for the reduction. This allows us to
describe the reduced space by means of Dirac's technique (see \cite{Dir}).
Then we show that different bialgebra structures can be used for the
reduction. The most important particular case corresponds to the bialgebra
structure related to the new Drinfeld realization of affine quantum groups.
It is this description of reduction that is important for quantization.

\subsection{ Drinfeld--Sokolov reduction for Poisson Lie groups}

\label{DS}

Recall the construction of the Drinfeld--Sokolov reduction for Poisson Lie
groups.

Let $G$ be a connected simply connected finite-dimensional complex
semisimple Lie group, ${\frak g}$ its Lie algebra. Fix a Cartan subalgebra $%
{\frak h}\subset {\frak g}\ $and let $\Delta $ be the set of roots of $%
\left( {\frak g},{\frak h}\right) .$ Choose an ordering in the root system;
let $\Delta _{+}$ be the system of positive roots, $\{\alpha _1,...,\alpha
_l\},$ $l=rank\,{\frak g},$ the set of simple roots and $H_1,\ldots ,H_l$
the set of simple root generators of ${\frak h}$.

Denote by $a_{ij}$ the corresponding Cartan matrix. Let $d_1,\ldots , d_l$
be coprime positive integers such that the matrix $b_{ij}=d_ia_{ij}$ is
symmetric. There exists a unique non--degenerate invariant symmetric
bilinear form $\left( ,\right) $ on ${\frak g}$ such that $(H_i ,
H_j)=d_j^{-1}a_{ij}$. It induces an isomorphism of vector spaces ${\frak h}%
\simeq {\frak h}^*$ under which $\alpha_i \in {\frak h}^*$ corresponds to $%
d_iH_i \in {\frak h}$. The induced bilinear form on ${\frak h}^*$ is given
by $(\alpha_i , \alpha_j)=b_{ij}$.

Let ${\frak b}$ be the positive Borel subalgebra and $\overline{{\frak b}}$
the opposite Borel subalgebra; let ${\frak n}=[{\frak b},{\frak b}]$ and $%
\overline{{\frak n}}=[\overline{{\frak b}},\overline{{\frak b}}]$ be their
nil-radicals. Let $H=\exp {\frak h},N=\exp {\frak n}, \overline{N}=\exp 
\overline{{\frak n}},B=HN,\overline{B}=H\overline{N}$ be the Cartan
subgroup, the maximal unipotent subgroups and the Borel subgroups of $G$
which correspond to the Lie subalgebras ${\frak h},{\frak n},\overline{%
{\frak n}},{\frak b}$ and $\overline{{\frak b}},$ respectively.

Let ${\bf g}=L{\frak g}$ be the loop algebra; we equip it with the standard
invariant bilinear form,

\begin{equation}  \label{prod}
\left\langle X,Y\right\rangle =Res_{z=0}\left( X\left( z\right) ,Y\left(
z\right) \right) dz/z.
\end{equation}

Define ${\bf b}=L{\frak b} ~,~\overline{{\bf b}}=L\overline{{\frak b}}~,~%
{\bf n}=Ln~,~\overline{{\bf n}}=L\overline{{\frak n}}~,~{\bf h}=L{\frak h}$
... The quotient algebras ${\bf b}/{\bf n}$, $\overline{{\bf b}}/\overline{%
{\bf n}}$ may be canonically identified with ${\bf h.}$

Denote by ${\cal G}~,~{\cal B}~,~\overline{{\cal B}}~,~{\cal N}~,~\overline{%
{\cal N}}~,~{\cal H}$ the corresponding loop groups. The groups $G~,~B~,~%
\overline{B}~,~ N~,~\overline{N}~,~H$ will be identified with the subgroups
of constant loops.

Introduce an operator $^{\theta }r\in {\rm End}~{\bf g}$ by 
\begin{equation}
^{\theta }r=P_{\overline{{\bf n}}}-P_{{\bf n}}+r^{0}~,~r^{0}={\frac{1+\theta 
}{1-\theta }}P_{{\bf h}},  \label{rmatrh}
\end{equation}
where $P_{\overline{{\bf n}}},~P_{{\bf n}}$ and $P_{{\bf h}}$ are the 
projection operators  onto ${\overline{{\bf n}}},~{\bf h}$ and ${{\bf h}}$, respectively, in the direct sum
$$
{\bf g}={\bf n}+\overline{{\bf n}}+{\bf h},
$$
and $\theta \in End~{\bf h}$ is a unitary
automorphism with respect to the scalar product (\ref{prod}) such that $det(\theta -1)\neq 0$. The operator (\ref{rmatrh})
satisfies the Classical Yang--Baxter equation.
Therefore, the Lie algebra ${\bf g}$ is equipped with the structure of a
factorizable Lie bialgebra. Moreover, in the notation of Theorem \ref{rmatrtheta}, we have ${\frak b}_{-}={\bf b}%
,{\frak b}_{+}=\overline{{\bf b}}$, ${\frak n}_{-}={\bf n},{\frak n}%
_{+}=\overline{{\bf n}}$ and $\theta_r=\theta$.

Fix $p\in {\Bbb C}$ , $\left| p\right| <1,$ and let $D_p $ be the
automorphism of ${\cal G}$ defined by $(D_pg)(z)=g(pz).$ We shall denote the
corresponding automorphism of the loop algebra ${\bf g}$ by the same letter.

Note that $D_{p}$ preserves the scalar product (\ref{prod}). Assume also
that $\theta $ commutes with $D_{p}$. Then the automorphism $D_{p}$ and the
r--matrix (\ref{rmatrh}) satisfy conditions (\ref{compat}). Therefore, we
can endow the group ${\cal G}$ with Poisson bracket (\ref{tau}) for $%
r=~^{\theta }r,~\sigma =D_{p}$. We denote the corresponding Poisson manifold
by ${\cal G}_{p,\theta }$. This manifold is equipped with the Poisson group
action (\ref{act}) of the Poisson Lie group ${\cal G}$:

\begin{equation}  \label{actn}
g\circ L = (g^{D_p})^{-1}Lg.
\end{equation}

Let $W$ be the Weyl group of $\left( {\frak g},{\frak h}\right) ;$ we shall
denote a representative of $w\in W$ in $G$ by the same letter. We also denote
$w(g)=wgw^{-1}$ for any $g\in G$. Let $%
s_1,...,s_l$ be the reflections which correspond to the simple roots; let $%
s=s_1s_2\cdot \cdot \cdot s_l$ be a Coxeter element.

The Drinfeld--Sokolov reduction for Poisson Lie groups is a Poisson
reduction for the Poisson manifold ${\cal G}_{p,\theta }$, where $\theta
=D_{p}\cdot s$, with respect to the gauge action (\ref{actn}) of the
unipotent group ${\cal N}\subset {\cal G}$. The result of the reduction is a
Poisson submanifold ${\cal S}\subset {\cal G}_{p,D_{p}\cdot s}/{\cal N}$
which has a nice geometric description.

Denote $N^{\prime }=\{v\in N;svs^{-1}\in \overline{N}\}$ , ${\cal N}^{\prime
}=\{v\in {\cal N};svs^{-1}\in \overline{{\cal N}}\}$ , $M^s=Ns^{-1}N$ , $%
{\cal M}^s={\cal N}s^{-1}{\cal N}$.

\begin{theorem}
{\bf (\cite{SS})}

(i) The gauge action (\ref{actn}) of the group ${\cal N}$ leaves the cell $%
{\cal M}^{s}\subset {\cal G}_{p,\theta }$ invariant. The action of ${\cal N}$
on ${\cal M}^{s}$ is free and ${\cal S}={\cal N}^{\prime }s^{-1}$ is a
cross-section of this action.

(ii) The subgroup ${\cal N}$ is admissible in ${\cal G}$ and hence ${\cal N}$%
--invariant functions form a Poisson subalgebra in the Poisson algebra $%
C^{\infty }({\cal G}_{p,\theta })$.

(iii) The quotient space ${\cal M}^{s}/{\cal N}$ =${\cal S}$ is a Poisson
submanifold in ${\cal G}_{p,\theta }/{\cal N}$ if and only if the
endomorphism $\theta $ is given by $\theta =s\cdot D_{p},$ where $s\in W$ is
a Coxeter element.
\end{theorem}

To simplify the notation we shall denote the Poisson manifold ${\cal G}_{p
,D_p \cdot s}$ by ${\cal G}_p $. The r--matrix $^{s\cdot D_p }r$ that enters 
the definition of the Poisson
structure of this manifold has the form:
\begin{equation}  \label{rs}
^{s\cdot D_p } r= P_{\overline{{\bf n}}}-P_{{\bf n}}+r^0 ,~~ r^0={\frac{%
1+s\cdot D_p }{1-s\cdot D_p}}P_{{\bf h}}.
\end{equation}

From (\ref{g*}) it follows that the
factorization problem (\ref{twfact}) for the manifold ${\cal G}_p ^{\prime}$
amounts to the relations :

\begin{eqnarray}  \label{factgp}
L=L^{D_p} _+ L_-^{-1}, \mbox{ where } L\in {\cal G}_p ^{\prime},~ (L_+,L_-)
\in {\cal G}^* , \\
L_\pm=h_\pm n_\pm ,~ h_\pm \in {\cal H},~h_-= s(h^{D_p} _+),~ n_+ \in 
\overline{{\cal N}} ,~ n_- \in {\cal N}.  \nonumber
\end{eqnarray}

We shall also consider the case of constant loops. The restriction of
operator (\ref{rmatrh}) to constant loops defines an r--matrix $^\theta r\in
End~{\frak g}$, where $\theta \in End~{\frak h}$. We denote the
corresponding Poisson--Lie group by $G$. Let $G_\theta$ be the manifold $G$
equipped with Poisson bracket (\ref{tau}), where $r=^\theta r,~\sigma =id$.
Then the action of $G$ on $G_\theta$ by conjugations is a Poisson group
action. We can formulate the finite--dimensional version of the previous
theorem.

\begin{corollary}
\label{cor1}

(i) The action of the group ${N}$ on $G_{\theta }$ by conjugations leaves
the cell ${M}^{s}\subset {G_{\theta }}$ invariant. The action of ${N}$ on ${M%
}^{s}$ is free and ${S}={N}^{\prime }s^{-1}$ is a cross-section of this
action.

(ii) The subgroup $N$ is admissible in $G$ and hence $N$--invariant
functions form a Poisson subalgebra in the Poisson algebra $C^{\infty }({G}%
_{\theta })$.

(iii) The quotient space ${M}^{s}/{N}$ =${S}$ is a Poisson submanifold in ${%
G_{\theta }}/{N}$ if and only if the endomorphism $\theta $ is given by $%
\theta =s$ , where $s\in W$ is the Coxeter element.
\end{corollary}

From the results of \cite{SS} (see section 4, Theorem 4.10) it follows that the 
Poisson structure of the reduced Poisson manifold $S$ may be described using 
a Poisson surjection $m: H\rightarrow S$ called the generalized Miura transform. 
For $\sigma
~=~id$ formula 4.3 in \cite{SS} implies that the Poisson structure of the 
manifold $H$ is trivial. Therefore the Poisson bracket on $S$ equals to zero
identically.

Note also that $N^{\prime }\subset N$ is an abelian subgroup, $\dim N^{s}=l$
(see \cite{BO}). So we arrive to
\begin{proposition}
The algebra of regular functions on the reduced space ${M}^{s}/{N}={N}%
^{\prime }s^{-1}$ is a polynomial algebra with $l$ generators; the reduced
Poisson structure is trivial.
\end{proposition}

\subsection{The structure of the r--matrix}

\label{srtheta}

In this section we show that r--matrix (\ref{rs}) is elliptic. Namely , we
shall express the kernel of its ``Cartan'' component $r^{0}$ by means of
theta functions.

The kernel of operator $r^{0}$ is given by\thinspace a formal power series

\begin{equation}
r^0(\frac zw)=\sum_{n=-\infty}^{\infty}{\frac{1+p^ns }{1-p^ns}}(\frac zw)^n ,
\end{equation}

where $s$ regarded as an element of ${\frak h}\otimes {\frak h}\simeq End~%
{\frak h}$ . 

Let $h$ be the Coxeter number of ${\frak g}$. Since $s^{h}=id,~$ $~r^{0}(z)$
may be represented in the following form:

\begin{equation}  \label{rzero}
r^0(z)=\sum_{m=0}^{h-1}s^m \sum_{n=-\infty}^{\infty}{\frac{p^{nm} }{1-p^{nh}}%
}z^n.
\end{equation}

According to Theorem 2.5 \cite{SS}, this formal power series satisfies the
following functional equation

\begin{equation}
\begin{array}{l}
\label{feq} r^0(z)= sr^0(pz)+s\delta (pz)+\delta (z), \mbox{ where } \delta
(z)=\sum_{n=-\infty}^\infty z^n .
\end{array}
\end{equation}

By successive application of this identity to $r^0(z)$ we obtain

\begin{equation}  \label{rd}
r^0(z)= r^0(p^hz)+ \delta (z)+\delta (p^hz)+2\sum_{k=1}^{h-1}s^k\delta
(p^kz).
\end{equation}

Therefore, if the formal power series (\ref{rzero}) has a nontrivial domain
of convergence it will define an $End~{\frak h}$ -valued elliptic function.
The delta functions in (\ref{rd}) indicate singularities of this function.
We shall show that this is indeed the case by comparing the series (\ref
{rzero}) with the Fourier expansion of an elliptic function. Recall that
every such function may be expressed via theta functions.

The standard theta function $\theta_t$ is defined by the relation

\begin{eqnarray}
\theta_t(u)= c\prod_{n=-\infty}^{\infty} (1-t^{2n-1}e^{2\pi
iu})(1-t^{2n-1}e^{-2\pi iu}) , \\
c=\prod_{n=1}^{\infty}(1-t^{2n}), ~t=e^{\pi i \xi} ,~ Im \xi >0  \nonumber
\end{eqnarray}

and satisfies the functional equations:

\begin{equation}
\begin{array}{l}
\theta_t(u+1)=\theta_t(u), \\ 
\\ 
\theta_t(u+\xi )=-t^{-1}e^{-2\pi iu}\theta_t(u).
\end{array}
\end{equation}

The only zeroes of $\theta_t$ are located at the points $m+(n+\frac 12)\xi ,
m,n \in {\Bbb Z}$.

In the stripe $|Im u |\leq \frac 12 Im \xi ,~ u\not{\in }\{m \pm \frac 12
\xi ,~ m \in {\Bbb Z}\}$ the ratio ${\frac{\theta_t^{\prime}}{\theta_t}}(u)$
has the following Fourier expansion (see \cite{BE}):

\begin{equation}  \label{func}
{\frac{\theta_t^{\prime}}{\theta_t}}(u)= {\frac{2\pi }{i}}%
\sum_{n=-\infty}^{\infty}{\frac{t^n }{1-t^{2n}}}e^{2\pi inu} .
\end{equation}

Put $z=e^{2\pi iu}$.

\begin{proposition}
For $u\in {\Bbb R}\backslash {\Bbb Z}$ the formal power series (\ref{rzero})
coincides with the Fourier expansion of the following elliptic function :

\begin{eqnarray*}
r^{0}(u)={\frac{i}{2\pi }}(2\sum_{m=1}^{h-1}s^{m}{\frac{\theta _{t}^{\prime }%
}{\theta _{t}}}(u+\xi (\frac{m}{h}-\frac{1}{2}))+ \\
{\frac{\theta _{t}^{\prime }}{\theta _{t}}}(u-\frac{\xi }{2})+{\frac{\theta
_{t}^{\prime }}{\theta _{t}}}(u+\frac{\xi }{2})),
\end{eqnarray*}

where $p$ and $t=e^{\pi i\xi }$ are related by $p^{\frac{h}{2}}=t$.
\end{proposition}

{\em Proof} follows from formula (\ref{func}) and expression (\ref{rzero})
for $r^0$.

\subsection{Dual pairs for some admissible actions and the Drinfeld--Sokolov
reduction}

Now we shall describe the Drinfeld--Sokolov reduction in terms of
constraints. For group actions admitting a moment map the constraints
imposed during reduction are given by matrix coefficients of the moment map.
The gauge action (\ref{actn}) of the subgroup ${\cal N}$ is not a Poisson
group action and hence it does not have a moment map in the usual sense.
Nevertheless, it is possible to define an analogue of the moment map in this
more general situation.

Let $A\times M \rightarrow M$ be a right Poisson group action of a
Poisson--Lie group $A$ on a Poisson manifold $M$. A subgroup $K\subset A$ is
called {\em admissible} if the set $C^\infty \left( M\right) ^K$ of $K$%
-invariants is a Poisson subalgebra in $C^\infty \left( M\right) .$

\begin{proposition}
\label{admiss}{\bf (\cite{RIMS}, Theorem 6 )}

Let $\left( {\frak a},{\frak a}^{*}\right) $ be the tangent Lie bialgebra of 
$A.$ A connected Lie subgroup $K\subset A$ with Lie algebra ${\frak k}%
\subset {\frak a}$ is admissible if ${\frak k}^{\perp }\subset {\frak a}^{*}$
is a Lie subalgebra.
\end{proposition}

Let $A\times M \rightarrow M$ be a right Poisson group action of a Poisson--Lie group $A$ on a manifold $M$.
Suppose that this action possesses a moment mapping $\mu : M\rightarrow A^*$.
Let $K$ be an admissible subgroup in $A$. Denote by $\frak k$ the Lie algebra of $K$. 
Assume that ${\frak k}^\perp \subset {\frak a}^*$ is a Lie subalgebra in ${\frak a}^*$. 
Suppose also that there is a splitting ${\frak a}^*={\frak t}\oplus {\frak k}^\perp$, and that 
$\frak t$ is a Lie subalgebra in ${\frak a}^*$. Then the linear space ${\frak k}^*$ is naturally
identified with $\frak t$.
Assume that $A^*$ is the semidirect 
product of the Lie subgroups $K^\perp , T$ corresponding to the Lie algebras
${\frak k}^\perp , {\frak t}$ respectively. Suppose that $K^\perp$ is a connected subgroup in $A^*$.
Fix the decomposition
$A^*=K^\perp T$ and denote by $\pi_{K^\perp} , \pi_{T}$ the projections onto
$K^\perp$ and $T$ in this decomposition. 
\begin{theorem}\label{QPmoment}
Define a map $\overline{\mu}:M\rightarrow T$ by
$$
\overline{\mu}=\pi_{T}\mu.
$$ 
Then
 
(i) 
$\overline{\mu}^*\left( C^\infty \left( T\right)\right)$ is a Poisson subalgebra in $C^\infty \left( M\right)$,
and hence one can equip $T$ with a Poisson structure such that $\overline{\mu}:M\rightarrow T$ is 
a Poisson map.

(ii)Moreover, the algebra $C^\infty \left( M\right) ^K$ is the centralizer of 
$\overline{\mu}^*\left( C^\infty \left( T\right)\right)$ in the Poisson algebra $C^\infty \left( M\right)$.
In particular, if $M/K$ is a smooth manifold the maps
\begin{equation}\label{dp}
\begin{array}{ccccc}
&  & M &  &  \\ 
& \stackrel{\pi }{\swarrow } &  & \stackrel{\overline{\mu}}{\searrow } &  ,\\ 
M/K &  &  &  & T
\end{array}
\end{equation}
form a dual pair.
\end{theorem}
{\em Proof.} (i)First, by Theorem 4.9 in \cite{Lu} there exists a Poisson bracket on $A^*$ such that 
$\mu :M\rightarrow A^*$ is a Poisson map. Moreover, we can choose this bracket to be the sum of the standard
Poisson--Lie bracket of $A^*$ and of a left invariant bivector on $A^*$. 
Denote by $A^*_M$ the manifold $A^*$ equipped with this Poisson structure.
Now observe that $T$ is identified with the quotient $K^\perp \setminus A^*_M$, where $K^\perp$ acts on
$A^*_M$ by multiplications from the left. Therefore to prove part (i) of the proposition it suffices to show
that $K^\perp$--invariant functions on $A^*_M$ form a Poisson subalgebra in $C^\infty(A^*_M)$.

Observe that since $A^*$ is a Poisson--Lie group and the Poisson structure of $A^*_M$ is obtained from that of
$A^*$ by adding a left--invariant term, the action of $A^*$ on $A^*_M$ by multiplications from the left is
a Poisson group action. Note also that $K^\perp$ is a connected subgroup in $A^*$ and $({\frak k}^\perp)^\perp 
\cong {\frak k}$ is a Lie subalgebra in $\frak a$. Therefore by Proposition \ref{admiss} $K^\perp$ is 
an admissible subgroup in $A^*$. Therefore 
$K^\perp$--invariant functions on $A^*_M$ form a Poisson subalgebra in $C^\infty(A^*_M)$, and hence
$\overline{\mu}^*\left( C^\infty \left( T\right)\right)$ is a Poisson subalgebra in $C^\infty \left( M\right)$.
This proves part (i).

(ii)By the definition of the moment map we have:
\begin{equation}\label{X5}
L_{\widehat X} \varphi =\langle \mu^*(\theta_{A^*}) , X \rangle (\xi_\varphi ) ,
\end{equation}
where $X \in {\frak a} , \widehat X$ is the corresponding vector field on $M$ and
$\xi_\varphi $ is the Hamiltonian vector field of $\varphi \in C^\infty (M)$. Since $A^*$ is the semidirect product of $K^\perp$ and $T$ the pullback of the right--invariant Maurer--Cartan form $\mu^*(\theta_{A^*})$ may be represented as follows:
$$
\mu^*(\theta_{A^*})= {\rm Ad}(\pi_{K^\perp}\mu )({\overline{\mu}}^*\theta_{T})+(\pi_{K^\perp}\mu )^*\theta_{K^\perp},
$$
where ${\rm Ad}(\pi_{K^\perp}\mu )({\overline{\mu}}^*\theta_{T})\in {\frak t},~(\pi_{K^\perp}\mu )^*\theta_{K^\perp}\in {\frak k}^\perp$.

Now let $X \in {\frak k}$. Then $\langle (\pi_{K^\perp}\mu )^*\theta_{K^\perp}),X\rangle =0$ and formula (\ref{X5}) takes the form:
\begin{equation}\label{+}
\begin{array}{l}
L_{\widehat X} \varphi =
\langle {\rm Ad}(\pi_{K^\perp}\mu )({\overline{\mu}}^*\theta_{T}),X \rangle (\xi_\varphi )=\\
\\
\langle {\rm Ad}(\pi_{K^\perp}\mu )(\theta_{T}),X \rangle ({\overline{\mu}}_*(\xi_\varphi )) .
\end{array}
\end{equation}

Since ${\rm Ad}(\pi_{K^\perp}\mu )$ is a non--degenerate transformation, $L_{\widehat X} \varphi =0$ for every $X\in {\frak k}$
if and only if ${\overline{\mu}}_*(\xi_\varphi )=0$, i.e. a function $\varphi \in C^\infty (M)$ is 
$K$--invariant if and only if $\{ \varphi ,\overline{\mu}^*(\psi) \}=0$ for every 
$\psi \in C^\infty (T)$. This completes the proof.

Recall that the notion of dual pairs serves to describe Poisson submanifolds
in the quotient $M/K$. If $M$ is symplectic, connected components of the
sets $\pi (\overline{\mu }^{-1}(x)),x\in T$ are symplectic leaves in $%
M/K $. This allows us to give an alternative description of the
Drinfeld--Sokolov reduction.

Using action (\ref{actn}) of the group ${\cal B}$ on ${\cal G}_{p}^{\prime } 
$ in the setting of the previous theorem we shall construct a dual pair for
the gauge action of the subgroup ${\cal N}\subset {\cal B}$. Observe that
according to part (iii) of Proposition \ref{bpm} $({\bf b},\overline{{\bf b}}%
)$ is a sub-bialgebra of $({\bf g},{\bf g}^{*})$. Therefore, ${\cal B}$ is a
Poisson Lie subgroup in ${\cal G}$. The restriction of the ${\cal B}$--gauge
action (\ref{actn}) to ${\cal G}_{p}^{\prime }$ has a moment map given by
Proposition \ref{gprime} and formula (\ref{factgp}) for the factorization
problem :

\begin{equation}
\mu_{{\cal B}}(L)=L_+ , \mbox{ where }L=L^{D_p} _+ L_-^{-1}, (L_+,L_-) \in 
{\cal G}^*.
\end{equation}

The orthogonal complement of the ${\bf n\subset b}$ in the dual space $%
\overline{{\bf b}}$ coincides with the ${\bf h}\subset \overline{{\bf b}}$.
Hence by Proposition \ref{admiss} ${\cal N}$ is an admissible subgroup in
the Poisson Lie group ${\cal B}$. Moreover, the dual group $\overline{{\cal B%
}}$ is the semidirect product of the Lie groups ${\cal H}$ and $\overline{%
{\cal N}}$ corresponding to the Lie algebras ${\bf n}^{\perp }={\bf h}$ and $%
\overline{{\bf n}}\cong {\bf n}^{*}$ , respectively. We conclude that all the
conditions of Theorem \ref{QPmoment} are satisfied with $A={\cal B},K={\cal N%
},A^{*}=\overline{{\cal B}},T=\overline{{\cal N}},K^{\perp }={\cal H}$.
It follows that Poisson manifold ${\cal G}_{p}^{\prime }$ possesses a dual
pair formed by the canonical projection onto the quotient ${\cal G}%
_{p}^{\prime }/{\cal N}$ and the map $\mu _{{\cal N}}$ ,

\begin{eqnarray}  \label{mun}
\mu_{{\cal N}}(L)=n_+, \mbox{ where } L=L^{D_p} _+ L_-^{-1}, (L_+,L_-) \in 
{\cal G}^* , \\
L_+=h_+n_+ ,~ h_+ \in {\cal H},~ n_+ \in \overline{{\cal N}}.  \nonumber
\end{eqnarray}

We shall describe the reduction with the help of the map $\mu_{{\cal N}}$.

Let $w_0\in W$ be the longest element; let $\tau \in Aut$ $\Delta _{+}$ be
the automorphism defined by $\tau \left( \alpha \right) =-w_0\cdot \alpha
,\alpha \in \Delta _{+}.$ Let $N_i\subset N$ be the 1-parameter subgroup
generated by the root vector corresponding to the root ${\tau \left( \alpha _i\right) }.$
Choose an element $u_i\in N_i,u_i\neq 1.$ Then we have (see \cite{st}) $%
w_0u_iw_0^{-1}\in Bs_iB.$ We may fix $u_i$ in such a way that $%
w_0u_iw_0^{-1}\in Ns_iN.$ Set $x=u_lu_{l-1}...u_1;$ then $f:=w_0xw_0^{-1}\in
Ns^{-1}N\cap \overline N.$

\begin{proposition}
\label{constrt} The set $\mu _{{\cal N}}^{-1}(f)$ is an open subset in $%
{\cal M}^{s}\cap {\cal G}_{p}^{\prime }$.
\end{proposition}

{\em Proof.} First , the space ${\cal M}^s$ is invariant with respect to the
following action of ${\cal H}$:

\[
h\circ L= hLs(h)^{-1}. 
\]

Indeed, let $L=vs^{-1}u ;~~ v,u \in {\cal N}$ be an element of ${\cal M}^s$.
Then

\begin{equation}
h\circ L=hvh^{-1}hs^{-1}s(h)^{-1}s(h)us(h)^{-1}=hvh^{-1}s^{-1}s(h)us(h)^{-1}.
\end{equation}

The r.h.s. of the last equality belongs to ${\cal M}^{s},$ because ${\cal H}$
normalizes ${\cal N}$.

Using definition (\ref{mun}) of $\mu _{{\cal N}}$ and formula (\ref{factgp})
for the factorization problem we can describe the level surface $\mu _{{\cal %
N}}^{-1}(f)$ as follows:

\[
\mu_{{\cal N}}^{-1}(f)= \{ h^{D_p}_+ fn_-^{-1} s(h^{D_p} _+ )^{-1}=
h^{D_p}_+ \circ (fn_-^{-1}) | n_- \in {\cal N} , h_+ \in {\cal H} \}. 
\]

Hence $\mu _{{\cal N}}^{-1}(f)\subset {\cal M}^{s}$. The dimension count
shows that $\mu _{{\cal N}}^{-1}(f)$ is open in ${\cal M}^{s}$.

This concludes the proof.

By Theorem \ref{QPmoment}, the image of the level surface $\mu _{{\cal N}%
}^{-1}(f)$ under the canonical projection ${\cal G}_{p}^{\prime }\rightarrow 
{\cal G}_{p}^{\prime }/{\cal N}$ is a Poisson submanifold in the quotient $%
{\cal G}_{p}^{\prime }/{\cal N}$. For the reasons which will be explained
later we denote this manifold by ${\cal G}_{p}^{\prime }/(F,\chi _{f})$.
Proposition \ref{constrt} implies that this manifold is open in ${\cal M}%
^{s}/{\cal N}$. Actually one can show that ${\cal G}_{p}^{\prime }/(F,\chi
_{f})$ is an open dense subset in ${\cal M}^{s}/{\cal N}$.

Now we want to describe the reduced manifold ${\cal G}_{p}^{\prime }/(F,\chi
_{f})$ in terms of constraints. To define the constraints we need the notion
of matrix coefficients of the group $\overline{{\cal N}}$ which is the
target space of the map $\mu _{{\cal N}}$.

Observe that one can define regular functions on $\overline{{\cal N}}$.
Indeed, let $\varphi $ be a regular function on $\overline{N}$. It may be
viewed as a function $\varphi ^{\prime }:\overline{{\cal N}}\rightarrow 
{\Bbb C}((z))$. For if $L=L(z)\in \overline{{\cal N}}$ then $\varphi
^{\prime }(L)=\varphi (L(z))$. The coefficients of the Laurent series $%
\varphi (L(z))$ are well defined regular functions on $\overline{{\cal N}}$.
In particular, one can define matrix coefficients of $\overline{{\cal N}}$.

It turns out that the constraints for the reduction , i.e. the matrix
coefficients of $\mu _{{\cal N}}$, are of the first class.

\begin{theorem}
\label{constrt1}

(i) Matrix coefficients of the map $\mu _{{\cal N}}$ form a Poisson
subalgebra $F$ in the Poisson algebra $C^{\infty }({\cal G}_{p}^{\prime })$.

(ii)Define a map $\chi _{f}:F\rightarrow {\Bbb C}$ by $\chi _{f}(\mu _{{\cal %
N}})=f$, where $\mu _{{\cal N}}$ should be viewed as a matrix of regular
functions on ${\cal G}_{p}^{\prime }$, the map $\chi _{f}$ is applied to the
matrix coefficients of $\mu _{{\cal N}}$. Then $\chi _{f}$ is a character of
the Poisson algebra $F$.
\end{theorem}

{\em Proof.}

(i)By part (i) of Theorem \ref{QPmoment} the algebra 
$\mu _{{\cal N}}^*\left( C^\infty \left( \overline{\cal N}\right)\right)$ is 
a Poisson subalgebra in $C^{\infty }({\cal G}_{p}^{\prime })$. We have to verify
that that the pullbacks of the matrix coefficients of $\overline{{\cal N}}$
form a Poisson subalgebra in $C^{\infty }({\cal G}_{p}^{\prime })$.

We shall use gradients of a function $\varphi\in C^\infty\left( {\cal G}%
_p ^{\prime}\right)$ with respect to the ${\cal G}^{*}$ group structure on
the set ${\cal G}_p ^{\prime}$ : 
\begin{eqnarray}
\langle X, \nabla \varphi (L_{+},L_{-})\rangle = \left( \frac d{ds}\right)
_{s=0}\varphi (e^{sX_{+}}L_{+},e^{sX_{-}}L_{-}),  \nonumber \\
\langle X, \nabla^{\prime} \varphi (L_{+},L_{-}) \rangle = \left( \frac
d{ds}\right) _{s=0}\varphi (L_{+}e^{sX_{+}},L_{-}e^{sX_{-}}),X \in {\bf g.}
\end{eqnarray}

We define $Z_\varphi \in {\bf g}$ by the following relation: 
\[
r_{+}Z_\varphi -D_p ^{-1}\cdot r_{-}Z_\varphi =\nabla^{\prime} \varphi. 
\]

Then the Poisson bracket on the Poisson submanifold ${\cal G}_p ^{\prime}$
may be represented as follows:

\begin{eqnarray}  \label{Pbracket}
\left\{ \varphi ,\psi \right\} _{D_p} (L_+,L_-)= \left\langle AdL_{+}\cdot
D_p ^{-1}\cdot Z_\varphi -AdL_{-}Z_\varphi ,\nabla \psi \right\rangle- \\
-\left\langle \nabla \varphi ,AdL_{+}\cdot D_p ^{-1}\cdot Z_\psi
-AdL_{-}Z_\psi \right\rangle .  \nonumber
\end{eqnarray}

This expression leads to the following explicit formula for the Poisson
bracket of functions which only depends on $n_{+}$ (see (\ref{mun})):

\begin{eqnarray}  \label{n+}
\left\{ \varphi ,\psi \right\} _{D_p} (L_+,L_-) = \left\langle{\frac{%
r_-^0+D_p r_+^0 }{r_-^0-D_p r_+^0}} P_{{\bf h}} (Adn_+^{-1} \nabla \varphi
), P_{{\bf h}} (Adn_+^{-1} \nabla \psi )\right\rangle + \\
+ \left\langle P_{{\bf b}} (Adn_+^{-1} \nabla \varphi ), (Adn_+^{-1} \nabla
\psi )\right\rangle - \left\langle (Adn_+^{-1} \nabla \varphi ), P_{{\bf b}}
(Adn_+^{-1} \nabla \psi )\right\rangle ,  \nonumber \\
\mbox{where }\left\langle \nabla \varphi (n_+) , X \right\rangle = \left(
\frac d{ds}\right) _{s=0}\varphi (e^{sX}n_+) ,\mbox{ for every } X \in 
\overline{{\bf n}}.  \nonumber
\end{eqnarray}

Clearly, such functions form a Poisson subalgebra.

Let $\varphi $ be a matrix coefficient of the map $\mu _{{\cal N}}$. Then
for every $X\in \overline{{\bf n}}$ $\left\langle \nabla \varphi
(n_{+}),X\right\rangle $ is an element of $F$ since $\left( \frac{d}{ds}%
\right) _{s=0}\varphi (e^{sX}n_{+})$ is a linear combination of matrix
coefficients. Therefore, for $\varphi ,\psi \in F$ the r.h.s. of (\ref{n+})
is an element of $F$.

(ii) We have to show that the Poisson bracket (\ref{n+}) vanishes when
restricted to the surface $\mu_{{\cal N}}^{-1}(f)$.

Recall that for every $X \in {\bf n}$ the action of the corresponding vector
field $\widehat X$ is given by formula (\ref{+}):

\begin{equation}
\begin{array}{l}
L_{\widehat X} \varphi= \langle Ad(\pi_{K^\perp}\mu )(\theta_{K^*}),X
\rangle ({\overline{\mu}}_*(\xi_\varphi)) .
\end{array}
\end{equation}

Now let $\varphi =\varphi (n_{+})$. The surface $\mu _{{\cal N}}^{-1}(f)$ is
stable , at least locally , under the gauge action of ${\cal N}$. This
implies $L_{\widehat{X}}\varphi |_{\mu _{{\cal N}}^{-1}(f)}=0$. Therefore, ${%
\overline{\mu }}_{*}(\xi _{\varphi })|_{\mu _{{\cal N}}^{-1}(f)}=0$, i.e.
the Poisson brackets of every function $\varphi =\varphi (n_{+})$ with the
matrix coefficients of $\mu _{{\cal N}}$ vanish on the constraint surface.
This completes the proof.

\begin{remark}
\label{rconstr}

According to Theorem \ref{constrt1}, the pair $(F,\chi _{f})$ is a system of
first class constraints for the Poisson manifold ${\cal G}_{p}^{\prime }$.
Using Dirac's technique (see\cite{Dir}), we can define the reduced manifold $%
{\cal G}_{p}^{\prime }/(F,\chi _{f})$ with the help of these constraints.
\end{remark}

\subsection{Deformation of the Poisson structure}

In this section we show that it is possible to use different bialgebra
structures to perform the Drinfeld--Sokolov reduction. First we prove that
the Poisson manifolds ${\cal G}_{p,\theta }$ are isomorphic for different $%
\theta $. For each $\theta $ we describe the reduced space ${\cal G}%
_{p}^{\prime }/(F,\chi _{f})$ using the corresponding Poisson manifold $%
{\cal G}_{p,\theta }^{\prime }$. Finally, we consider an important case
associated with the bialgebra structure related to the new Drinfeld
realization of affine quantum groups.

Let $\theta ,\theta ^{\prime }\in End~{\bf h}$ be two unitary endomorphisms
such that $det(\theta -1)\neq 0,det(\theta ^{\prime }-1)\neq 0$. Suppose
that they commute with $D_{p}$ and with each other. According to the results
of Section \ref{DS}, every such endomorphism defines a factorizable
bialgebra structure on ${\bf g}$. As a consequence, we obtain two Poisson
manifolds ${\cal G}_{p,\theta },{\cal G}_{p,\theta ^{\prime }}$ equipped
with Poisson group ${\cal G}$--actions (\ref{actn}). The restrictions of
these actions to the submanifolds ${\cal G}_{p,\theta }^{\prime },{\cal G}%
_{p,\theta ^{\prime }}^{\prime }$ possess moment maps (see Proposition \ref
{gprime}). The results of \cite{A} imply that the Poisson manifolds ${\cal G}%
_{p,\theta }^{\prime },{\cal G}_{p,\theta ^{\prime }}^{\prime }$ are
isomorphic. Moreover, the isomorphism is given by the gauge action (\ref
{actn}), where $g$ depends on $L$. We can define $g(L)$ more precisely.

Let $L$ be an element of ${\cal G}_{p,\theta }^{\prime }$. From (\ref{imbd})
it follows that the factorization (\ref{twfact}) for the manifold ${\cal G}%
_{p,\theta }^{\prime }$ amounts to the relations :

\begin{equation}  \label{facth}
\begin{array}{l}
L=L_+^{D_p} L_-^{-1}~,\mbox{ where } L_\pm =h_\pm n_\pm , \\ 
\\ 
h_\pm \in {\cal H} ~,~n_+ \in \overline{{\cal N}} ~,~ n_- \in {\cal N} ~,~
h_\pm =e^{^{\theta }r_\pm^0X} ~,~ X \in {\bf h}.
\end{array}
\end{equation}

\begin{proposition}
Let $A\in End~{\bf h}$ be an endomorphism commuting with $\theta ,\theta
^{\prime }$ and $D_{p}$. The map

\begin{equation}
{\cal G}_{p,\theta }^{\prime }\rightarrow {\cal G}_{p,\theta ^{\prime
}}^{\prime }:L\mapsto t^{D_{p}}Lt^{-1}=L^{\prime }~,~t=e^{AX},  \label{iso}
\end{equation}
where $X$ is given by (\ref{facth}), is an isomorphism of the Poisson
manifolds if and only if $A$ satisfies the equation:

\begin{equation}
AA^{*}{\frac{D_{p}-1}{^{\theta }r_{-}^{0}-D_{p}~^{\theta }r_{+}^{0}}}%
+A-A^{*}=~^{\theta ^{\prime }}r_{+}^{0}-~^{\theta }r_{+}^{0}.
\end{equation}

Let $L^{\prime }={L_{+}^{\prime }}^{D_{p}}{L_{-}^{\prime }}^{-1}~,~L_{\pm
}^{\prime }=h_{\pm }^{\prime }n_{\pm }^{\prime }~,~h_{\pm }^{\prime }\in 
{\cal H}~,~n_{+}^{\prime }\in \overline{{\cal N}}~,~n_{-}^{\prime }\in {\cal %
N}$ be the factorization (\ref{facth}) of $L^{\prime }$ as an element of the
manifold ${\cal G}_{p,\theta ^{\prime }}^{\prime }$. Then in terms of the
components $h_{\pm }^{\prime },~n_{\pm }^{\prime }$ the map (\ref{iso}) has
the form:

\[
h_{\pm }^{\prime }=e^{^{\theta ^{\prime }}r_{\pm }^{0}Y},~~~Y={\frac{%
D_{p}A-A+D_{p}~^{\theta }r_{+}^{0}-~^{\theta }r_{-}^{0}}{D_{p}~^{\theta
^{\prime }}r_{+}^{0}-~^{\theta ^{\prime }}r_{-}^{0}}}X,
\]

\[
n_{+}^{\prime }=e^{KX}n_{+}e^{-KX},~~~n_{-}^{\prime
}=e^{D_{p}KX}n_{-}e^{-D_{p}KX},
\]

\[
K\in End~{\bf h},~~~K={\frac{A+~^{\theta }r_{+}^{0}-~^{\theta ^{\prime
}}r_{+}^{0}}{D_{p}~^{\theta ^{\prime }}r_{+}^{0}-~^{\theta ^{\prime
}}r_{-}^{0}}},
\]
where $n_{\pm }$ and $X$ are given by (\ref{facth}).

The operator $K$ satisfies the equation:

\[
(1-D_{p})KK^{*}+D_{p}K-K^{*}={\frac{^{\theta }r_{+}^{0}-~^{\theta ^{\prime
}}r_{+}^{0}}{^{\theta ^{\prime }}r_{+}^{0}-{D_{p}}^{-1}~{^{\theta ^{\prime
}}r_{-}^{0}}}}.
\]
\end{proposition}

{\em Proof} is provided by direct calculation using formula (\ref{Pbracket})
for the Poisson bracket on ${\cal G}_{p,\theta} ^{\prime}$.

The main idea of this proposition is that one can use different Poisson
structures for the Drinfeld--Sokolov reduction. From Theorem \ref{constrt1}
, remark \ref{rconstr} and the previous proposition for $\theta
^{\prime}=D_p \cdot s$ we obtain the following

\begin{proposition}
\label{defred} Let $L\in {\cal G}_{p,\theta }^{\prime
},L=L_{+}^{D_{p}}L_{-}^{-1}$ be an element of ${\cal G}_{p,\theta }^{\prime }
$. Fix an operator $A$ defined in the previous theorem. Consider the map

\[
\mu _{{\cal N}}^{\theta ,K}:{\cal G}_{p,\theta }^{\prime }\rightarrow 
\overline{{\cal N}};~~\mu _{{\cal N}}^{\theta ,K}(L)=e^{KX}n_{+}e^{-KX}
\]

which is the composition of the isomorphism (\ref{iso}) and the moment map $%
\mu _{{\cal N}}$. Then

(i) Matrix coefficients of $\mu _{{\cal N}}^{\theta ,K}$ form a Poisson
subalgebra $F^{\theta ,K}$ in $C^{\infty }({\cal G}_{p,\theta }^{\prime })$.

(ii) The map $\chi _{f}^{\theta ,K}:F^{\theta ,K}\rightarrow {\Bbb C}$
defined by $\chi _{f}^{\theta ,K}(\mu _{{\cal N}}^{\theta ,K})=f$ is a
character of the Poisson algebra $F^{\theta ,K}$. Thus $(F^{\theta ,K},\chi
_{f}^{\theta ,K})$ is a system of the first class constraints.

(iii) The reduced Poisson manifold ${\cal G}_{p,\theta }^{\prime
}/(F^{\theta ,K},\chi _{f}^{\theta ,K})$ is isomorphic to ${\cal G}%
_{p}^{\prime }/(F,\chi _{f})$.
\end{proposition}

Now we shall describe the reduced space ${\cal G}_{p}^{\prime }/(F,\chi
_{f}) $ using the bialgebra structure related to the 'new Drinfeld
realization' of affine quantum groups \cite{nr}. Recall that this bialgebra
structure is factorizable, the corresponding r--matrix is given by:

\begin{equation}  \label{rdm}
^Dr=P_{\overline{{\bf n}}} -P_{{\bf n}} +^Dr^0 ,
\end{equation}

where

\[
^Dr^0=P^0_+ -P^0_-, 
\]

and $P^0_\pm$ are the projection operators onto $z{\frak h}[[z]]$ and $z^{-1}%
{\frak h}[[z^{-1}]]$ , respectively , in the direct sum

\begin{equation}  \label{hsum}
{\bf h}=z^{-1}{\frak h}[[z^{-1}]]\stackrel{\cdot}{+}{\frak h}\stackrel{\cdot%
}{+}z^{-1}{\frak h}[[z^{-1}]] .
\end{equation}

Denote by ${\cal G}_{p,D}$ the manifold ${\cal G}$ equipped with Poisson
bracket (\ref{tau}), where $r=~^{D}r,~\sigma =D_{p}$. The previous
proposition cannot be directly applied to the manifold ${\cal G}_{p,D}$,
since the r--matrix $^{D}r$ is not of form (\ref{rmatrh}). But it may be
obtained by a limit procedure from an r--matrix of this type.

Indeed , consider the unitary automorphism $\theta \in End ~{\frak h}$
defined by $( \theta h)(z)=h(uz), u \in {\Bbb C} , |u|<1$. Then the kernel
of the ``Cartan'' component $^{\theta }{r^0}$ of the corresponding r--matrix (%
\ref{rmatrh}) is

\begin{equation}  \label{lim}
^{\theta }{r^0}(z)=\sum_{n=-\infty}^{\infty}t{\frac{1+u^n }{1-u^n}}z^n,
\end{equation}

where $t\in {\frak h}\otimes {\frak h}$ is the Casimir element of ${\frak h}$%
... Consider the limit $u\rightarrow 0$. The only part of the r--matrix
depending on $u$ is $^{\theta }{r^{0}}$. From formula (\ref{lim}) it follows
that the limit of $^{\theta }{r^{0}}$ is a well--defined skew--symmetric
operator on ${\bf h}$ which coincides with $^{D}r^{0}$. Therefore, the
r--matrix $^{\theta }r$ degenerates into $^{D}r$.

By the limit procedure from Proposition \ref{defred} we get

\begin{theorem}
\label{TK} Let $L\in {\cal G}_{p,D}^{\prime }$ be an element of the manifold 
${\cal G}_{p,D}^{\prime }$ factorized as in (\ref{twfact}) : $%
L=L_{+}^{D_{p}}L_{-}^{-1}~,~L_{\pm }=h_{\pm }n_{\pm }~,~h_{\pm }\in {\cal H}%
~,~n_{+}\in \overline{{\cal N}}~,~n_{-}\in {\cal N}~,~h_{\pm }=e^{^{D}r_{\pm
}^{0}X}~,~X\in {\bf h}$ .

Define the map $\mu _{{\cal N}}^{D,K}:{\cal G}_{p,D}^{\prime }\rightarrow 
\overline{{\cal N}}$ by

\begin{equation}
\mu _{{\cal N}}^{D,K}(L)=e^{KX}n_{+}e^{-KX},  \label{cal}
\end{equation}

where $K\in End~{\bf h}$ is an endomorphism commuting with $s,~D_{p}$ and
satisfying the equation:

\begin{equation}
(1-D_{p})KK^{*}+D_{p}K-K^{*}=-{\frac{D_{p}s}{1-s}}P_{+}^{0}-{\frac{1}{1-s}}%
P_{-}^{0}+\frac{1}{2}{\frac{1+s}{1-s}}P_{0},  \label{K}
\end{equation}

$P_{0}$ is the projection operator onto ${\frak h}$ in the direct sum (\ref
{hsum}).

Then

(i) Matrix coefficients of $\mu _{{\cal N}}^{D,K}$ form a Poisson subalgebra 
$F^{D,K}$ in $C^{\infty }({\cal G}_{p,D}^{\prime })$.

(ii) The map $\chi _{f}^{D,K}:F^{D,K}\rightarrow {\Bbb C}$ defined by $\chi
_{f}^{D,K}(\mu _{{\cal N}}^{D,K})=f$ is a character of the Poisson algebra $%
F^{D,K}$. Thus $(F^{D,K},\chi _{f}^{D,K})$ is a system of the first class
constraints.

(iii) The reduced Poisson manifold ${\cal G}_{p,D}^{\prime }/(F^{D,K},\chi
_{f}^{D,K})$ is isomorphic to ${\cal G}_{p}^{\prime }/(F,\chi _{f})$.
\end{theorem}

Note that the equation (\ref{K}) does not contain elliptic functions in the
r.h.s..

Similarly to this theorem, we obtain the following description of the
reduction for constant loops (see corollary \ref{cor1}).

\begin{corollary}
\label{cor2} Denote by $G_{D}^{\prime }$ the manifold $G$ equipped with
Poisson bracket (\ref{tau}), where $r$ is the restriction of the r--matrix $%
^{D}r$ to constant loops and $\sigma ~=~id$. Let $L\in G_{D}^{\prime }$ be
an element of the manifold $G_{D}^{\prime }$ factorized as in (\ref{twfact})
: $L=L_{+}L_{-}^{-1}~,~L_{\pm }=h_{\pm }n_{\pm }~,~h_{\pm }\in {H}%
~,~n_{+}\in \overline{N}~,~n_{-}\in {N}~,~h_{\pm }=e^{^{D}r_{\pm
}^{0}X}~,~X\in {h}$ .

Define the map $\mu _{N}^{D,K}:{G}_{D}^{\prime }\rightarrow \overline{N}$ by

\begin{equation}
\mu _{N}^{D,K}(L)=e^{KX}n_{+}e^{-KX},  \label{mufin}
\end{equation}

where $K\in End~{h}$ is an endomorphism commuting with $s$ and satisfying
the equation:

\begin{equation}
K-K^{*}=\frac{1}{2}{\frac{1+s}{1-s}}.  \label{K0}
\end{equation}

Then

(i) Matrix coefficients of $\mu _{N}^{D,K}$ form a Poisson subalgebra $%
F^{0,K}$ in $C^{\infty }({G}_{D}^{\prime })$.

(ii) The map $\chi _{f}^{0,K}:F^{0,K}\rightarrow {\Bbb C}$ defined by $\chi
_{f}^{0,K}(\mu _{N}^{D,K})=f$ is a character of the Poisson algebra $F^{0,K}$%
... Thus $(F^{0,K},\chi _{f}^{0,K})$ is a system of the first class
constraints.

(iii) The reduced Poisson manifold ${G}_{D}^{\prime }/(F^{0,K},\chi
_{f}^{0,K})$ is isomorphic to an open Poisson submanifold in $S$.
\end{corollary}

\section*{Discussion}

In conclusion we briefly discuss quantization of the Poisson--Lie version of
the Drinfeld--Sokolov reduction.

First observe that the affine quantum group $U_{q}(\widehat{{\frak g}})$
with central charge $q^{c}=p$ is a quantization of the Poisson manifold $%
{\cal G}_{p,D}^{\prime }$ (see \cite{JF},\cite{ext}). The algebra $U_{q}(%
\widehat{{\frak g}})$ contains the quantum group $U_{q}({\frak g})$ which is
a quantization of the Poisson manifold ${\ G}_{D}^{\prime }$ (see \cite{FRT}%
, \S 2; \cite{dual}, \S 3).

For simplicity we shall consider the reduction for constant loops (see
corollary \ref{cor2}) in detail. Let us examine map (\ref{mufin}). The group 
$\overline{N}$ is unipotent and may be identified with its Lie algebra $%
\overline{{\frak n}}$ by means of the exponential map. Let $f_{i},i=1,\ldots
,l$ be the simple root generators of $\overline{{\frak n}}$. An element $%
n_{+}\in \overline{N}$ may be written as $e^{\sum_{i=1}^{l}f_{i}\varphi
_{i}+\psi },$ where $\varphi _{i}\in {\Bbb C}$ and $\psi $ is a term of
higher order with respect to the principal grading of $\overline{{\frak n}}$.

Using the notation of corollary \ref{cor2} , we can expand $X\in {\frak h}$
with respect to the basis of root generators as follows: $%
X=\sum_{i=1}^{l}H_{i}\psi _{i}$. Consider $\varphi _{i},\psi _{i},i=1,\ldots
,l$ as functions $\varphi _{i},\psi _{i}:\overline{B}\rightarrow {\Bbb C}%
,\varphi _{i}(L_{+})=\varphi _{i},\psi _{i}(L_{+})=\psi _{i}$. Then the map (%
\ref{mufin}) induces the following mapping of functions:

\begin{equation}  \label{anz}
(\mu^{D,K}_{N})^*(\varphi_i)=e^{\sum_{j=1}^l\langle \alpha_i , KH_j\rangle
\psi_j}\varphi_i .
\end{equation}
The r.h.s. of (\ref{anz}) is an element of the algebra of constraints.

Functions $\varphi _{i},\psi _{i}$ correspond to simple positive root
generators and the fundamental weight generators of $U_{q}({\frak g})$,
respectively. In \cite{S1} we have shown that transformation (\ref{anz}) has
an exact quantum counterpart (see formula (19) in \cite{S1}). Quantum
constraints defined in \cite{S1} form a subalgebra in $U_{q}({\frak g})$.
This subalgebra possesses a character which is a quantum counterpart of $%
\chi _{f}^{0,K}$. Equation (\ref{K0}) appears in the quantum case as well
(equation (18) in \cite{S1}). This allows us to establish an exact
correspondence between the classical and the quantum pictures. A detailed
analysis of the quantum reduction is contained in \cite{thes}, Chapter 4.

A similar situation is observed in the affine case. Matrix coefficients of
the components $X,n_{\pm }$ introduced in Theorem \ref{TK} correspond to the
loop generators of $U_{q}(\widehat{{\frak g}})$ in Drinfeld's 'new
realization' (see \cite{JF},\cite{nr}). Quantum analogues of map (\ref{cal})
, equation (\ref{K}) and the character $\chi _{f}^{D,K}$ are defined in
terms of Drinfeld's 'new realization' as well (see formula (23), equation
(33) and discussion after Proposition 8 in \cite{S1}).

Finally, we remark that a quantization of the reduced space ${\cal G}%
_{p,D}^{\prime }/(F^{D,K},\chi _{f}^{D,K})$ may be obtained by applying the
homological reduction procedure proposed in \cite{S2} for arbitrary systems
of the first--class constraints. This program will be realized in a
subsequent paper.

\end{document}